\documentclass[10pt]{article}
\usepackage{latexsym,amssymb}
\def\demo{\noindent{\bf Proof. }}
\def\QED{\hfill$\Box$}

\newtheorem{Theorem}{Theorem}[section]
\newtheorem{Lemma}[Theorem]{Lemma}
\newtheorem{Corollary}[Theorem]{Corollary}
\newtheorem{Proposition}[Theorem]{Proposition}
\newtheorem{Remark}[Theorem]{Remark}
\newtheorem{Example}[Theorem]{Example}

\begin{document}

\topmargin3mm
\hoffset=-1cm
\voffset=-1.5cm
\begin{center}

{\large\bf Normalization of monomial ideals and Hilbert functions}\\

\vspace{6mm}


\addtocounter{footnote}{1}


\medskip

Rafael
H. Villarreal\footnote{Partially supported
by CONACyT grant 49251-F and SNI, M\'exico.}\\

{\small Departamento de
Matem\'aticas}\vspace{-1mm}\\
{\small Centro de Investigaci\'on y de Estudios
Avanzados del
IPN}\vspace{-1mm}\\
{\small Apartado Postal
14--740}\vspace{-1mm}\\
{\small 07000 Mexico City, D.F.}\vspace{-1mm}\\
{\small
e-mail: {\tt
vila@math.cinvestav.mx}}\vspace{4mm}

\end{center}

\date{}

\begin{abstract}
We study the normalization of a monomial ideal, and show how to compute
its Hilbert function (using Ehrhart polynomials) if the ideal is zero
dimensional. A positive 
lower bound for the second coefficient of the Hilbert polynomial is
shown.
\end{abstract}

\section{Introduction}

Let $R=k[x_1,\ldots,x_d]$ be a polynomial ring over a field $k$ and
let $I$  be a monomial ideal of $R$ minimally generated by 
$x^{v_1},\ldots,x^{v_q}$. As usual for 
$a=(a_i)$ in $\mathbb{N}^d$ we set $x^a=x_1^{a_1}\cdots x_d^{a_d}$.  
If ${\mathcal R}$ is the Rees algebra of $I$, ${\mathcal R} =
 R[It]$, we call its integral closure $\overline{{\mathcal R}}$ the
 {\em normalization} of $I$. 
This algebra has for components the integral closures of the powers
of $I$:
\[ {\mathcal R} = R\oplus It\oplus\cdots\oplus I^it^i\oplus\cdots 
\subset R\oplus \overline{I}t\oplus\cdots\oplus
\overline{I^i}t^i\oplus\cdots=\overline{{\mathcal R}}.\]

In our situation ${\mathcal R}\subset\overline{\mathcal R}$ is a
finite extension. By a result of Vasconcelos \cite[Theorem~7.58]{bookthree}, 
the $I$-filtration ${\mathcal
F}=\{\overline{I^i}\}_{i=0}^\infty$ stabilizes 
for $i\geq d$, i.e., 
$\overline{I^i}=I\, \overline{I^{i-1}}$ for $i\geq d$. 
We complement this result by showing that if $\deg(x^{v_i})=r$ for
all $i$, then $\mathcal F$ stabilizes for 
$i$ greater or equal than the minimum of ${\rm rank}(v_1,\ldots,v_q)$
and  $d-\lfloor d/r\rfloor+1$ (Proposition~\ref{aug6-05} 
and Corollary~\ref{dec14-06}).  

If $\dim(R/I)=0$, we  are interested in studying the 
{\it Hilbert function\/} of $\mathcal F$:
$$
f(n)=\ell_R(R/\overline{I^n})=c_dn^d+c_{d-1}n^{d-1}+\cdots+c_1n+c_0
\ \ (c_i\in\mathbb{Q};\ n\gg 0). 
$$

We will
express $f(n)$ as a difference of two Ehrhart polynomials
(Proposition~\ref{may29-1-05}) and show
a lower bound for $c_{d-1}$  (Proposition~\ref{boletin}). In
particular we obtain an efficient way of computing the Hilbert
function of $\mathcal F$ using integer programming methods. As
an application we show that $e_0(d-1)-2e_1\geq d-1$, where $e_i$ is the
$i${\it th} Hilbert coefficient of $f$. For monomial
ideals this improves the inequality $e_0(d-1)\geq 2e_1$ given by Polini, 
Ulrich and Vasconcelos \cite[Theorem~3.2]{PUV} that holds for an arbitrary
$\mathfrak{m}$-primary ideal $I$ of a regular local ring 
$(R,\mathfrak{m})$. This inequality turns out to be useful to bound 
the length of divisorial chains for classes of Rees algebras 
\cite[Corollary~3.5]{PUV}.

In the sequel we use \cite{BHer,Schr} as references for 
standard terminology and notation on commutative algebra and 
polyhedral geometry. We denote the
set of non-negative real (resp. integer, rational) numbers  by $\mathbb{R}_+$ 
(resp. $\mathbb{N}$, $\mathbb{Q}_+$).

\section{Normalization of monomial ideals}\label{normalization}

To avoid repetition, we continue to use the notation and definitions
used in the introduction.  

\begin{Proposition}\label{aug6-05} Let $r_0$ be the rank of the matrix
$(v_1,\ldots,v_q)$. If $v_1,\ldots,v_q$ lie in 
a hyperplane of $\mathbb{R}^d$ not containing the 
origin, then $\overline{I^b}=I\overline{I^{b-1}}$ for $b\geq r_0$.  
\end{Proposition}

\demo  Let 
$\mathbb{Q}_+{\mathcal A}'$ be the cone in $\mathbb{Q}^{d+1}$ generated by the set
$$
{\mathcal A}'=\{(v_1,1),\ldots,(v_q,1),e_1,\ldots,e_d\},
$$
where $e_i$ is the $i${\it th} unit vector in 
$\mathbb{Q}^{d+1}$. Assume $b\geq r_0$. Notice that we invariably have 
$I\overline{I^{b-1}}\subset \overline{I^b}$. To show the reverse
inclusion take $x^\alpha\in\overline{I^b}$, i.e., 
$x^{m\alpha}\in I^{bm}$ for some $0\neq m\in\mathbb{N}$. Hence 
$(\alpha,b)\in\mathbb{Q}_+{\mathcal A}'$. Applying Carath\'eodory's  
theorem for cones \cite[Corollary 7.1i]{Schr}, we can write 
$$
(\alpha,b)=\lambda_1(v_{i_1},1)+\cdots+\lambda_r(v_{i_r},1)
+\mu_1e_{j_1}+\cdots+\mu_se_{j_s}\ \ \ 
(\lambda_\ell,\mu_k\in\mathbb{Q}_+),
$$
where $\{(v_{i_1},1),\ldots,(v_{i_r},1),e_{j_1},\ldots,e_{j_s}\}$ 
is a linearly independent set contained in ${\mathcal A}'$. 
Notice that $v_{i_1},\ldots,v_{i_r}$ are also linearly independent 
because they lie in a hyperplane not containing the origin. 
Hence $r\leq r_0$. Since $b=\lambda_1+\cdots+\lambda_r$, we obtain 
that $\lambda_\ell\geq 1$
for some $\ell$, say $\ell=1$. Then 
$$
\alpha=v_{i_1}+(\lambda_1-1)v_{i_1}+\lambda_2v_{i_2}+\cdots+
\lambda_rv_{i_r}+\mu_1e_{j_1}+\cdots+\mu_se_{j_s},
$$
and consequently $x^\alpha\in I\overline{I^{b-1}}$. \QED

\begin{Remark}\rm  In this proof we may replace
$\mathbb{Q}$ by $\mathbb{R}$ because 
according to \cite[p.~219]{monalg} we have the equality 
$\mathbb{Z}^{d+1}\cap \mathbb{R}_+{\mathcal A}'=
\mathbb{Z}^{d+1}\cap \mathbb{Q}_+{\mathcal A}'$.
\end{Remark}

\begin{Lemma}\label{irrrepecve} Let 
${\mathcal A}'=\{e_1,\ldots,e_d\}\cup\{(a_1,\ldots,a_d,1)
\vert\, a_i\in\mathbb{N};\ \sum_ia_i=r\}$, where 
$r\geq 2$ is an integer. Then the irreducible representation 
of the cone\/ $\mathbb{R}_+{\mathcal A}'$, as an intersection of
closed halfspaces, is given by 
$$
\mathbb{R}_+{\mathcal A}'=H_{e_1}^+\cap\cdots\cap H_{e_d}^+
\cap H_{e_{d+1}}^+\cap H_{a}^+,
$$
where $a=(1,\ldots,1,-r)$ and $H_a^+=
\{x\in\mathbb{R}^{d+1}\, 
\vert\, \langle x,a\rangle\geq 0\}$.
\end{Lemma}

\demo  We set ${\mathcal A}=
\{e_1,\ldots,e_d,re_1+e_{d+1},\ldots,re_d+e_{d+1}\}$ and
$N=\{e_1,\ldots,e_{d+1},a\}$. The cone $\mathbb{R}_+{\mathcal A}'$ has
dimension $d+1$ and one has the equality 
$\mathbb{R}_+{\mathcal A}'=\mathbb{R}_+{\mathcal A}$.  
Thus it suffices 
to prove that $F$ is a facet of $\mathbb{R}_+{\mathcal A}$ 
if and only if $F=H_b\cap{\mathbb R}_+{\mathcal A}$ for some $b\in N$. 
Let $1\leq i\leq d+1$. Consider the following sets of vectors:  
$$
\Gamma_i=\{e_1,e_2,\ldots,\widehat{e}_i,\ldots, e_d,e_{d+1}\}\ \mbox{ and }\ 
\Gamma=\{re_1+e_{d+1},\ldots,re_d+e_{d+1}\},
$$
where $\widehat{e}_i$ means to omit $e_i$ from the list. Since
$\Gamma_i$ and $\Gamma$ are linearly independent, we obtain that 
$F=H_b\cap{\mathbb R}_+{\mathcal A}$ is a facet for $b\in N$, i.e., 
$\dim(F)=d$ and $\mathbb{R}_+{\mathcal A}\subset H_b^+$. Conversely 
let $F$ be a facet of the cone 
${\mathbb R}_+{\mathcal A}$. There 
are linearly independent vectors $\alpha_1,\ldots ,\alpha_d\in{\mathcal A}$ 
 and  $0\neq b=(b_1,\ldots,b_{d+1})\in{\mathbb R}^{d+1}$ 
such that 
\begin{enumerate}
\item[\rm (i)\ \ ] $F={\mathbb R}_+{\mathcal A}\cap H_b$,
\item[\rm (ii)\ ] $\mathbb{R}\alpha_1+\cdots+\mathbb{R}\alpha_d=H_b$, 
and 
\item[\rm (iii)] ${\mathbb R}_+{\mathcal A}\subset H^+_b$.
\end{enumerate}
Since ${e}_{1},\ldots,{e}_{d}$ are in 
$\mathcal A$, by (iii) one has:
\begin{equation}\label{ago30-1}
\langle {e}_{i},b\rangle=b_i\geq 0\ 
\mbox{ for }\ i=1,\ldots,d.
\end{equation}
Set ${\mathcal B}=\{\alpha _1,\ldots , \alpha_d\}$ and consider 
the matrix $M$ whose rows are the vectors in $\mathcal B$. To finish
the proof we need only show that there exists $c\in N$ such that
$H_b=H_c$. Consider the following 
cases. Case (1): If the $i${\it th\/} 
column of $M$ is zero for 
some $1\leq i\leq d+1$, we set $c={e}_i$. Case (2): If ${\mathcal
B}=\{re_1+e_{d+1},\ldots,re_d+e_{d+1}\}$,  
then we set $c=a$. Case (3): Now we assume that  
$$
{\mathcal B}=\{{e}_{i_1}, \ldots ,{e}_{i_s}, re_{j_1}+e_{d+1},\ldots,
re_{j_t}+e_{d+1}\},
$$ 
where $s, t>0$, $s+t=d$, $1\leq i_1<\cdots< i_s\leq d$, 
$1\leq j_1<\cdots< j_t\leq d$, and $M$ has all its columns 
different from zero. Since $b_{i_1}=0$, using 
$$
\begin{array}{ccl}
\langle re_{i_1}+e_{d+1},b\rangle&=&b_{d+1}\geq 0,\\
\langle re_{j_1}+e_{d+1},b\rangle&=&rb_{j_1}+b_{d+1}=0,
\end{array}
$$
and Eq.(\ref{ago30-1}), we obtain $b_{d+1}=0$. 
Then $e_{d+1}\in H_b$. It follows 
readily that $H_b$ is generated, as a vector space, by the set 
$\{e_1,e_2,\ldots,\widehat{e}_i,\ldots, e_d,e_{d+1}\}$ 
for some $1\leq i\leq d$. Thus in this case we set $c=e_i$.
\QED 

\medskip

Let $r\geq 2$ be an integer. For the rest of this section we make two
assumptions: (i) $R[t]$ has the grading
$\delta$ induced by setting $\delta(x_i)=1$ and $\delta(t)=1-r$, and
(ii) $\deg(x^{v_i})=r$ for all $i$.  Thus $\mathcal{R}=R[It]$ becomes a standard
graded $k$-algebra. In this case $\mathcal{R}$ and
$\overline{\mathcal{R}}$ have rational Hilbert series. The degree as a
rational function of the Hilbert series of $\mathcal R$, denoted by
$a(\mathcal{R})$, is  
called the $a$-{\it invariant\/} of $\mathcal{R}$. 

For use below recall that the $r${\it th} {\it Veronese ideal\/} of
$R$ is the ideal generated 
by all monomials of $R$ of degree $r$. If $x^a$ is a monomial we 
set $\log(x^a)=a$. 

\begin{Proposition}\label{inv-veron} 
Let $J$ be the $r${\it th} Veronese 
ideal of $R$ and let $S=R[Jt]$ be its Rees algebra. {\rm (a)} If
$r\geq d$, then $a(S)=-2$.  
{\rm (b)} If $2\leq r<d$ and 
$d=qr+s$, 
where $0\leq s<r$, then 
$$
a(S)=\left\{
\begin{array}{lcl}
-(q+2)& \mbox{ if }& s\geq 2,\\ 
-(q+1)& \mbox{ if }&s = 0\mbox{ or }s=1. 
\end{array}\right.
$$
\end{Proposition}

\demo  Let ${\mathcal A}'$ be as in Lemma~\ref{irrrepecve}. As $S$ is normal, 
according to a formula of Danilov-Stanley \cite{BHer}, the
canonical module $\omega_{S}$ of $S$ can be expressed as
\begin{equation}\label{dan-stan}
\omega_{S}=(\{x^at^b\vert\,
(a,b)\in {\mathbb N}{\mathcal A}'\cap({\mathbb R}_+{\mathcal
A}')^{\rm o}\})=(\{x^at^b\vert\,
(a,b)\in \mathbb{Z}^{d+1}\cap({\mathbb R}_+{\mathcal
A}')^{\rm o}\}),
\end{equation} 
where $({\mathbb R}_+{\mathcal A}')^{\rm o}$ 
denotes the relative interior of ${\mathbb R}_+{\mathcal A}'$ and 
${\mathbb N}{\mathcal A}'$ is the subsemigroup of $\mathbb{N}^{d+1}$
generated by ${\mathcal A}'$. In our situation recall that $a(S)
=-\min\{i\, \vert\, (\omega_S)_i\neq 0\}$.

Let $m\in\omega_S$. We can write $m=x^a(x^bt^c)$,  
where $x^bt^c=(f_1t)\cdots(f_ct)$ 
and $f_i$ is a monomial of degree $r$ for all $i$. 
Notice that $\delta(m)=|a|+c$, where $a=(a_i)$ and
$|a|=a_1+\cdots+a_d$. Since $\log(m)=(a+b,c)$ is in  
the interior 
of the cone ${\mathbb R}_+{\mathcal A}'$, using 
Lemma~\ref{irrrepecve} one 
has $c\geq 1$, $a_i+b_i\geq 1$ for all $i$, and 
$|a|+|b|\geq rc+1$. As 
$|b|=rc$, altogether we get:
\begin{equation}\label{jun5}
|a|+|b|\geq d\mbox{ and }|a|\geq 1.
\end{equation}
In particular $\delta(m)\geq 2$. This shows the inequality $a(S)\leq -2$ 
because $m$ was an arbitrary monomial in $\omega_S$. To prove (a) notice
that by Lemma~\ref{irrrepecve} the monomial 
$m_1=x_1^{r-d+2}x_2\cdots x_dt$ is in $\omega_S$ and $\delta(m_1)=2$.
Hence $a(S)=-2$. To prove (b) there are three cases 
to consider. We only show the 
case $s\geq 2$, the cases $s=1$ and $s=0$ can be shown similarly.

Case $s\geq 2$: First we show that $\delta(m)\geq q+2$. If $c>q$, 
then from Eq.(\ref{jun5}) we get $\delta(m)\geq q+2$. Assume $c\leq
q$. One has the inequality: 
\begin{equation}\label{eqnum1}
r(q-c)+s\geq (q-c)+2.
\end{equation}
From Eq.(\ref{jun5}) one has $|a|+|b|=|a|+rc\geq d=rq+s$. 
Consequently
\begin{equation}\label{eqnum2}
\delta(m)=|a|+c\geq r(q-c)+s+c.
\end{equation}
Hence from Eqs.(\ref{eqnum1}) and (\ref{eqnum2})  we get $\delta(m)\geq q+2$. 
Therefore one has the inequality $a(S)\leq -(q+2)$, to show 
equality it suffices to prove that the monomial
$$
m_2=x_1^2x_2^2\cdots x_{r-s+1}^2 x_{r-s+2} \cdots x_d t^{q + 1}
$$
is in $\omega_S$ and has degree $q+2$. An easy calculation 
shows that $\delta(m_2)=q+2$. Finally let us see that $m_2$ is in $\omega_S$ via 
Lemma~\ref{irrrepecve}. That the entries of $\log(m_2)$ satisfy 
$X_i>0$ for all $i$ is clear.  The inequality
\begin{equation}\label{dec13-06}
X_1+X_2+\cdots+X_d > rX_{d+1}, 
\end{equation}
after making $X_i$ equal to the $i${\it th} entry of $\log(m_2)$,  
transforms into
$$
2(r-s+1)+(d-(r-s+1)) >r(q+1),
$$
but the left hand side is $r(q+1)+1$, hence ${\rm log}(m_2)$ 
satisfies Eq.~(\ref{dec13-06}). Hence ${\rm log}(m_2)$ is in the
interior of $\mathbb{R}_+{\mathcal A}'$, i.e., $m_2\in\omega_S$.\QED 

\medskip

The next result sharpen 
\cite[Theorem~3.3]{BVV} for the class 
of ideals generated by monomials of the same degree. 

\begin{Proposition}\label{mar2-02} If $2\leq r<d$, then the 
normalization $\overline{\mathcal R}$ 
of $I$ is generated as an ${\mathcal R}$-module by elements $g\in R[t]$ 
of $t$-degree at most $d-\lfloor d/r\rfloor$. 
\end{Proposition}

\demo  Set $f_i=x^{v_i}$ for $i=1,\ldots,q$. Consider the subsemigroup 
$C$ of $\mathbb{N}^{d+1}$ generated by the set $\{e_1,\ldots,e_d, 
(v_1,1),\ldots,(v_q,1)\}$ and the subgroup $\mathbb{Z}C$ generated by
$C$. Since ${\mathbb Z}C={\mathbb
Z}^{d+1}$, the normalization of $I$ can be expressed as:
\[
\overline{\mathcal R}=k[\{x^at^b\vert\, (a,b)\in{\mathbb Z}^{d+1}\cap{\mathbb
R}_+C\}].
\]
Let $m=x^at^b$ be a monomial of $\overline{\mathcal R}$ with
$(a,b)\neq 0$.  
We claim that $\delta(m)\geq b$. To show this inequality write
$$
(a,b)=\lambda_1e_1+\cdots+\lambda_de_d+
\mu_1\log(f_1t)+\cdots+\mu_q\log(f_qt),  
$$
where $\lambda_i\geq 0,\mu_j\geq 0$ for all $i,j$. Hence
$$
\begin{array}{ccl}
|a|=\lambda_1+\cdots+\lambda_d+(\mu_1+\cdots+\mu_q)r&\mbox{and}&
b=\mu_1+\cdots+\mu_q. 
\end{array}
$$
Consequently $\delta(m)=|a|+(1-r)b=
(\lambda_1+\cdots+\lambda_d)+b\geq b$. We may assume that $k$ is
infinite. There
is a Noether normalization
$A=k[z_1,\ldots,z_{d+1}]\stackrel{\varphi}{\hookrightarrow}{\mathcal 
R}$ such that $z_1,\ldots,z_{d+1}\in {\mathcal R}_1$. If $\psi$ is the
inclusion from ${\mathcal R}$ to $\overline{\mathcal R}$,
note that 
$A\stackrel{\psi\varphi}{\longrightarrow}\overline{\mathcal R}
$ is a Noether normalization. By 
\cite{Ho1}, the ring $\overline{\mathcal R}$ is Cohen-Macaulay.
Hence $\overline{\mathcal R}$ is a
free $A$-module that according to  \cite[Proposition~2.2.14]{monalg}
can be written as 
\begin{equation}\label{mar3-02}
\overline{\mathcal R}=Am_1\oplus\cdots\oplus A m_n,
\end{equation}
where $m_i=x^{\beta_i}t^{b_i}$. Set $h_i=|\{j\, \vert\,
\delta(m_j)=i\}|$. 
Using that the length is additive we
obtain the following expression for the Hilbert series of
$\overline{\mathcal{R}}$:
$$
H(\overline{\mathcal R},z)=
\sum_{i=0}^n\frac{z^{\delta(m_i)}}{(1-z)^{d+1}}=
\frac{h_0+h_1z+\cdots+h_sz^s}{(1-z)^{d+1}}.
$$
Recall that $a(\overline{\mathcal R})
=-\min\{i\, \vert\, (\omega_{\overline{\mathcal R}})_i\neq 0\}$, 
where $\omega_{\overline{\mathcal R}}$ is the canonical module of 
$\overline{\mathcal{R}}$. Let $J$ be the $r${\it th} Veronese 
ideal of $R$ and let $S=R[Jt]$ be its Rees algebra. Notice that
$\overline{\mathcal R}\subset S$ because $R[It]\subset S$ and $S$ 
is normal. Since $\dim(\overline{\mathcal R})=\dim(S)=d+1$, from the
Danilov-Stanley formula (see Eq.~(\ref{dan-stan})) it is seen that 
$a(\overline{\mathcal R})\leq a(S)$; see the proof of 
\cite[Proposition~3.5]{BVV}. Therefore using
Proposition~\ref{inv-veron} we get:  
$$a(\overline{\mathcal R})=s-(d+1)\leq a(R[Jt])\leq 
-\left\lfloor{d}/{r}\right\rfloor-1,$$
and $s\leq d-\lfloor d/r\rfloor$. Altogether if 
$m_i=x^{\beta_i}t^{b_i}$, one has:
\begin{equation}\label{dec15-06}
b_i\leq
 \delta(m_i)\leq s\leq d+1+a(R[Jt])\leq d-\lfloor d/r\rfloor.
\end{equation}
Therefore the $t$-degree of $m_i$ is less or equal than 
$d-\lfloor d/r\rfloor$, as required. \QED 

\begin{Proposition}\label{aug6-1-05} $\overline{I^{b}}=I\overline{I^{b-1}}$ 
for $b\geq d+2+a(R[Jt])$. 
\end{Proposition}

\demo  It suffices to prove the inclusion
$\overline{I^{b}}\subset I\overline{I^{b-1}}$. Let
$x^a\in\overline{I^b}$, i.e., $m=x^at^b\in\overline{\mathcal R}$. From  
Eq.~(\ref{mar3-02}) and noticing that $A\subset\mathcal R$, we can
write $m=(x^\gamma t^c)m_i$ for some $i$, where
$m_i=x^{\beta_i}t^{b_i}$ and $x^{\gamma}\in I^c$. Using
Eq.~(\ref{dec15-06}) gives $c\geq 1$. Thus $x^a\in
I^c\overline{I^{b_i}}$. To complete the proof notice that 
$I^c\overline{I^{b_i}}=I(I^{c-1}\overline{I^{b_i}})\subset
I\overline{I^{b_i+c-1}}=I\overline{I^{b-1}}$.  
\QED 

\begin{Corollary}\label{dec14-06}
$\overline{I^{b}}=I\overline{I^{b-1}}$ for $b\geq d-\lfloor
d/r\rfloor+1$. 
\end{Corollary}

\demo  By Proposition~\ref{inv-veron} one has $a(R[Jt])\leq 
-\lfloor d/r\rfloor-1$. Hence the result follows applying
Proposition~\ref{aug6-1-05}.  
\QED 

\section{Zero dimensional monomial ideals}

Let $R=k[x_1,\ldots,x_d]$ be a polynomial ring over a field $k$, 
with $d\geq 2$, and let $I$ be a zero dimensional monomial ideal
of $R$ minimally generated by $x^{v_1},\ldots,x^{v_q}$. Here we will
study the integral closure of  
the powers of $I$ and its Hilbert function.

We may assume 
that $v_i=a_ie_i$ for $1\leq i\leq d$, where $a_1,\ldots,a_d$ are
positive integers and $e_i$ is the $i${\it th\/} unit vector of 
$\mathbb{Q}^d$. Set $\alpha_0=(1/a_1,\ldots,1/a_d)$. We may also assume that 
$\{v_{d+1},\ldots,v_s\}$ is the set of $v_i$ such that 
$\langle v_i,\alpha_0\rangle <1$, and $\{v_{s+1},\ldots,v_{q}\}$ 
is the set of $v_i$ such that $i>d$ and $\langle
v_i,\alpha_0\rangle\geq 1$. Consider the convex polytopes in
$\mathbb{Q}^d$:
$$
P:={\rm conv}(v_1,\ldots,v_s),\ \ \ S:={\rm conv}(0,v_1,\ldots,v_d)
=\{x\, \vert\, x\geq 0;\, \langle x,\alpha_0\rangle\leq 1\},
$$
and the rational convex polyhedron $Q:=\mathbb{Q}_+^d+{\rm
conv}(v_1,\ldots,v_q)$. 

\begin{Proposition}\label{may29-05} $\overline{I^n}=(\{x^a\vert\,
a\in{nQ}\cap\mathbb{Z}^d\})$ for $0\neq n\in\mathbb{N}$.
\end{Proposition}

\demo  Let $x^\alpha\in\overline{I^n}$, i.e., 
$x^{m\alpha}\in I^{nm}$ for some $0\neq m\in\mathbb{N}$. Hence 
$$
\alpha/n\in{\rm conv}(v_1,\ldots,v_q)+\mathbb{Q}_+^d=Q
$$
and $\alpha\in nQ\cap\mathbb{Z}^d$. Conversely 
let $\alpha\in{nQ}\cap\mathbb{Z}^d$. It is seen that 
$x^{m\alpha}\in{I}^{nm}$ for some $0\neq
m\in\mathbb{N}$, this yields $x^\alpha\in\overline{I^n}$.\QED 

\medskip

Let us give a simpler expression for $Q$. From the equality  
$$
\mathbb{Q}_+^d+{\rm conv}(v_1,\ldots,v_d) 
=\{x\, \vert\, x\geq 0;\, \langle x,\alpha_0\rangle\geq 1\},
$$
we get that $v_i\in\mathbb{Q}_+^d+P$ for $i=1,\ldots,q$. Using
the finite basis theorem for polyhedra \cite[Corollary~7.1b]{Schr} we have that
$\mathbb{Q}_+^d+P$ is a convex set. Hence 
$Q\subset \mathbb{Q}_+^d+P$, and consequently we obtain the equality
\begin{equation}\label{dec16-2-06}
Q=\mathbb{Q}_+^d+P.
\end{equation}

\begin{Corollary} If $\langle v_i,\alpha_0\rangle\geq 1$ for all $i$,
then $\overline{I^n}=\overline{(x_1^{a_1},\ldots,x_d^{a_d})^n}$ 
for $n\geq 1$.
\end{Corollary}

\demo 
It follows at once from Proposition~\ref{may29-05} and
Eq.~(\ref{dec16-2-06}). Notice that in this case 
$P={\rm conv}(a_1e_1,\ldots,a_de_d)$. \QED  

\medskip

The {\it Hilbert function\/} of the filtration 
${\mathcal F}=\{\overline{I^n}\}_{n=0}^\infty$ is defined as
$$
f(n)=\ell(R/\overline{I^n})=\dim_k(R/\overline{I^n});\ \ \
n\in\mathbb{N}\setminus\{0\};\ \ f(0)=0.
$$
As usual $\ell(R/\overline{I^n})$ denotes the length of 
$R/\overline{I^n}$ as an $R$-module. 
For simplicity we call $f$ the {\it Hilbert function\/} of $I$. 

\begin{Corollary}\label{incredibles} 
$\ell(R/\overline{I^n})=|\mathbb{N}^d\setminus nQ|$ for $n\geq 1$.
\end{Corollary}

\demo  The length of $R/\overline{I^n}$ equals the dimension of
$R/\overline{I^n}$ as a $k$-vector space. By
Proposition~\ref{may29-05} the set 
${\mathcal B}=\{\overline{x^c}\vert\, c\notin nQ\}$ is precisely the set of
standard monomials of $R/\overline{I^n}$. Thus $\mathcal B$ is a
$k$-vector space basis of $R/\overline{I^n}$, and the equality
follows. \QED  

\medskip

The function $f$ is a polynomial function of degree $d$:
$$
f(n)=c_dn^d+c_{d-1}n^{d-1}+\cdots+c_1n+c_0\ \ \ \ (n\gg 0), 
$$
where $c_0,\ldots,c_d\in\mathbb{Q}$ and $c_d\neq 0$. The polynomial 
$c_dx^d+\cdots+c_0$ is called the {\it Hilbert polynomial\/} of $\mathcal
F$. One has the equality
$d!c_d=e(I)=e(\overline{I})$, where $e(I)$ is the multiplicity of
$I$, see \cite{multical}. We will
express $f(n)$ as a difference of two Ehrhart polynomials and then show
a positive lower bound for $c_{d-1}$.  

The {\it Ehrhart function\/} of $P$ is 
the numerical function
$\chi_P\colon\mathbb{N}\rightarrow\mathbb{N}$ given 
by $\chi_P(n)=|\mathbb{Z}^d\cap nP|$. This is a
polynomial function of degree $d_1=\dim(P)$:
$$
\chi_P(n)=b_{d_1}n^{d_1}+\cdots+b_1n+b_0\ \ \ (n\gg 0),
$$
where $b_i\in\mathbb{Q}$ for all $i$. The polynomial  
$E_P(x)=b_{d_1}x^{d_1}+\cdots+b_1x+b_0$ 
is called the {\it Ehrhart polynomial\/} of
$P$. 

\begin{Remark}\label{ehrhart-prop}\rm Some well known properties of $E_P$ are 
(see \cite{BHer}): 
\begin{enumerate}
\item $b_{d_1}={\rm vol}(P)$, 
where ${\rm vol}(P)$ denotes the relative volume of $P$.         
\item $b_{d_1-1}=(\sum_{i=1}^s{\rm vol}(F_i))/2$ where
$F_1,\ldots,F_s$ are the facets of $P$.   
\item $\chi_P(n)=E_P(n)$ for all integers $n\geq 0$. In 
particular $E_P(0)=1$. 
\item {\it Reciprocity law of Ehrhart\/}: $E_P^{\rm
o}(n)=(-1)^{d}E_P(-n)\ \ \forall\, n\geq 1,$\\  
where $E_P^{\rm o}(n)=|\mathbb{Z}^d\cap(nP)^{\rm o}|$ 
and $(nP)^{\rm o}$ is the relative interior of $nP$. 
\end{enumerate}
\end{Remark}

\begin{Lemma}\label{dec19-06} $P=S\cap Q$. 
\end{Lemma}

\demo  Clearly $P\subset S\cap Q$. Conversely let $z=(z_i)\in
S\cap Q$. Assume that $z\notin P$. By the separating hyperplane
theorem \cite[Theorem~3.23]{korte}, there are $0\neq b=(b_i)\in\mathbb{R}^d$ and
$c\in \mathbb{R}$ such that $\langle b,v_i\rangle\leq c$ for
$i=1,\ldots,s$ and $\langle b,z\rangle>c$. Assume $c>0$. Since 
$b_ia_i\leq c$ for all $i$ and $\langle
\alpha_0,z\rangle\leq 1$, we get $\langle b,z\rangle\leq \langle
\alpha_0,z\rangle c\leq c$, a contradiction. If $c=0$, then $b_i\leq
0$ for all $i$ and $\langle b,z\rangle\leq 0$, a contradiction. If 
$c<0$, we write $z=\delta+p$, for some $\delta\in\mathbb{Q}_+^d$ and
$p\in P$. Then $c<\langle b,z\rangle=\langle
b,\delta\rangle+\langle
b,p\rangle\leq \langle b,\delta\rangle+c$. Thus $0<\langle
b,\delta\rangle$, a contradiction because $b_i\leq 0$ for all $i$.
\QED 

\begin{Proposition}\label{may29-1-05} $f(n)=E_{S}(n)-E_P(n)$ for
$n\in\mathbb{N}$. In particular 
$$f(n)=c_dn^d+c_{d-1}n^{d-1}+\cdots+c_1n+c_0\mbox{ for }
n\in\mathbb{N} \mbox{ and }c_0=0.
$$
\end{Proposition}

\demo  Since $E_{P}(0)=E_{S}(0)=1$, we get the equality at $n=0$. 
Assume $n\geq 1$. Using Lemma~\ref{dec19-06}, we get the 
decomposition $Q=(\mathbb{Q}_+^d\setminus S)\cup{P}$. Hence
$$
nQ=(\mathbb{Q}_+^d\setminus nS)\cup{nP}\ \Longrightarrow\ 
\mathbb{N}^d\setminus
nQ=[\mathbb{N}^d\cap(nS)]\setminus[\mathbb{N}^d\cap(nP)].
$$
Therefore by Corollary~\ref{incredibles} we obtain 
$f(n)=E_{S}(n)-E_P(n)$. \QED

\begin{Example}\label{may19-05}\rm Let
$I=(x_1^4,x_2^5,x_3^6,x_1x_2x_3^2)$. Notice that 
$$
P={\rm conv}((4,0,0),(0,5,0),(0,0,6),(1,1,2)).
$$
Using 
{\it Normaliz\/} \cite{B}, to compute the Ehrhart polynomials of $S$ and
$P$,  we get 
\begin{eqnarray*}
&f(n)=E_{S}(n)-E_{P}(n)=(1+6n+19n^2+20n^3)\ \ \ \ \ \ \ \ \ \ 
\ \ \ \ \ \ \ \ \ \ \ \ \ \ \ \ \ \ \ \ \ \ \ \ \ &\\ 
&-(1+(1/6)n+(3/2)n^2+(13/3)n^3)=(35/6)n+(35/2)n^2+(47/3)n^3.&
\end{eqnarray*}
\end{Example}

\begin{Theorem}{\rm \cite[Theorem~7.58]{bookthree}}\label{wolmer-angra}
$\overline{I^b}=I\overline{I^{b-1}}$ for $b\geq d$.  
\end{Theorem}

\begin{Remark}\rm  We can use polynomial interpolation
together with Theorem~\ref{wolmer-angra} and
Proposition~\ref{may29-1-05} to determine
$c_1,\ldots,c_d$, see Example~\ref{polyinterpex}. 
\end{Remark}

\begin{Example}\label{polyinterpex}\rm Let
$I=(x_1^{10},x_2^8,x_3^5)$. 
Using {\it CoCoA\/} \cite{CNR} we obtain that the values of 
$f$ at $n=0,1,2,3$ are $0,112,704,2176$. By polynomial interpolation we
get:
$$
f(n)=\ell(R/\overline{I^n})=({200}/{3})n^3+40x^2+({16}/{3})n,\ \
\forall\, n\geq 0.
$$ 
\end{Example}

\begin{Lemma}\label{may20-05} Let $\alpha=(\alpha_i)$ and $\beta=(\beta_i)$ be 
two vectors in $\mathbb{Q}_+^d$ such that $\alpha_i=\beta_i$ for
$i=1,\ldots,d-1$, $\beta_d>\alpha_d$ and $\langle
\beta,\alpha_0\rangle<1$. Then
\begin{enumerate}
\item[\rm (a)] $\beta\in{\rm conv}(v_1,\ldots,v_d,\alpha)$.
\item[\rm (b)] If $\alpha_i>0$ for $i=1,\ldots,d-1$, then 
$\beta\in{\rm conv}(v_1,\ldots,v_d,\alpha)^{\rm o}$.
\item[\rm (c)] If $\alpha_i>0$ for $i=1,\ldots,d-1$ and 
$\alpha\in P$, then $\beta\in P^{\rm o}$. 
\end{enumerate}
\end{Lemma}

\demo  (a) To see that $\beta$ is a convex combination of 
$v_1,\ldots,v_d,\alpha$ we set:
\begin{eqnarray*}
s&=&\sum_{i=1}^d\alpha_i/a_i=\langle\alpha_0,\alpha\rangle<1,\ \ 
\mu=1-\left[\frac{\beta_d-\alpha_d}{a_d(1-s)}\right]> 0, 
\\
\lambda_i&=&(1-\mu){\alpha_i}/{a_i}\geq 0,\ \ \ i=1,\ldots,d-1,\\
\lambda_d&=&({\beta_d-\mu\alpha_d})/{a_d}=(({\beta_d-\alpha_d})/{a_d})+
{\alpha_d}(1-\mu)/{a_d}>0.
\end{eqnarray*}
Then $\beta=\lambda_1v_1+\cdots+\lambda_dv_d+\mu\alpha$ 
and $\lambda_1+\cdots+\lambda_d+\mu=1$, as required.

(b) Set $V=\{v_1,\ldots,v_d,\alpha\}$ and $\Delta={\rm conv}(V)$.
Since $V$ is affinely independent, $\Delta$ is a
$d$-simplex. From \cite[Theorem~7.3]{Bron}, the facets of 
$\Delta$ are precisely those sets of the form ${\rm conv}(W)$, where
$W$ is a subset of $V$ having $d$ points. If $\beta$ is not in the 
interior of $\Delta$, then $\beta$ must lie
in its boundary by (a). Therefore $\beta$ lies in some facet of
$\Delta$, which rapidly yields a contradiction.  

(c) By part (b) we get 
$\beta\in{\rm conv}(v_1,\ldots,v_d,\alpha)^{\rm o}\subset P^{\rm o}$,
as required. \QED 

\medskip

\noindent {\it Notation} The relative boundary of $P$ will be 
denoted by $\partial P$. 

\begin{Lemma}\label{may20-1-05} If $\alpha\in\partial P\setminus{\rm
conv}(v_1,\ldots,v_d)$ and $\alpha_i>0$ for $i=1,\ldots,d$, 
then the vector $\alpha'=(\alpha_1,\ldots,\alpha_{d-1},0)$ is not 
in $P$. 
\end{Lemma}

\demo  Notice that $\langle\alpha,\alpha_0\rangle<1$. If $\alpha'\in P$,
then by Lemma~\ref{may20-05}(c) we obtain $\alpha\in P^{\rm o}$, a contradiction.
Thus $\alpha'\notin P$.  \QED 

\medskip

For use below we set 
\begin{eqnarray*}
K_i&=&\{(a_i)\in S\vert\,
a_i=0\}={\rm conv}(\{v_1,\ldots,v_d,0\}\setminus\{v_i\});\ \
1\leq i\leq d,\\
H&=&{\rm conv}(v_1,\ldots,v_d);\ \ K=(\cup_{i=1}^d K_i)\setminus H; \
\ L=\partial P\setminus H,\ \mbox{ if
} H\subsetneq P. 
\end{eqnarray*}

Consider the map $\psi\colon L\rightarrow K$ given by
$$
\psi(\alpha)=\left\{
\begin{array}{ll}
\alpha,&\mbox{ if } \alpha_i=0 \mbox{ for some } 1\leq i\leq d,\\
(\alpha_1,\ldots,\alpha_{d-1},0),& \mbox{ if } \alpha_i>0\mbox{ for
all } 
1\leq i\leq d.
\end{array}
\right.
$$
Take $\alpha\in L$. Then $\langle \alpha,\alpha_0\rangle <1$. 
Since $\partial{P}\subset
P\subset S$ it is seen that $\psi(\alpha)\in K$. Indeed if
$\alpha_i=0$ for some $i$, then $\psi(\alpha)=\alpha\in K_i\setminus
H$. If $\alpha_i>0$ for all $i$, then $\alpha$ is a convex
combination of $v_1,\ldots,v_d,0$. Hence $\psi(\alpha)$ is a convex
combination of $v_1,\ldots,v_{d-1},0$ and $\psi(\alpha)\in
K_d\setminus H$.

\begin{Lemma}\label{injective} $\psi$ is injective.
\end{Lemma}

\demo  Let $\alpha=(\alpha_i),\beta=(\beta_i)\in L$. 
Assume $\psi(\alpha)=\psi(\beta)$. 
If $\alpha_i=0$ for some $i$ and $\beta_j=0$ for some $j$, then 
clearly $\alpha=\beta$. If $\beta_i>0$ for $i=1,\ldots,d$ and 
$\alpha_j=0$ for some $j$, then $\alpha_d=0$ and $\alpha_i=\beta_i$
for $i=1,\ldots,d-1$, by Lemma~\ref{may20-1-05} 
we can readily see that this case cannot
occur. If $\alpha_i\beta_i>0$ for all $i$, then 
$\alpha=\beta$ by Lemma~\ref{may20-05}(c). \QED 

\medskip

Let us introduce some more notation. We set
\begin{eqnarray*}
{\mathcal A}_i&=&\{v_j\vert\, 1\leq j\leq s;\, x_i\notin{\rm
supp}(x^{v_j})\};\\
P_i&=&{\rm conv}({\mathcal A}_i); \ \ 
H_i={\rm conv}(\{v_1,\ldots,v_d\}\setminus\{v_i\})\subset P_i\subset
K_i.
\end{eqnarray*}

\begin{Lemma}\label{may26-05} 
$\partial P\cap K_i=P_i$ for $i=1,\ldots,d$.
\end{Lemma}

\demo  For simplicity of notation assume $i=1$. Let $\alpha=(\alpha_i)
\in\partial P\cap K_1$, then $\alpha\in P$ and $\alpha_1=0$. Since
$\alpha$ is a convex combination of $v_1,\ldots,v_s$ it follows
rapidly that $\alpha$ is a convex combination of ${\mathcal A}_1$, i.e., 
$\alpha\in P_1$. Conversely let $\alpha\in P_1$. Clearly 
$\alpha\in K_1\cap P$ because ${\mathcal A}_1\subset K_1\cap P$. 
Assume that $\alpha\notin\partial P$. Then $\alpha\in P^{\rm o}$.  
If $\dim(P)=d-1$, we have that $P={\rm conv}(v_1,\ldots,v_d)$ and $P$
is a simplex. Thus by \cite[Theorem~7.3]{Bron}, the facets
of $P$ are $F_1,\ldots,F_d$, where  
$F_i={\rm conv}(v_1,\ldots,v_{i-1},v_{i+1},\ldots,v_d)$. The relative boundary
of $P$ is equal to $F_1\cup\cdots\cup F_d$. Hence
$\alpha\notin F_i$ for all $i$ and we can
write $\alpha=\lambda_1v_1+\cdots+\lambda_dv_d$, where 
$\sum_{i=1}^d\lambda_i=1$ and $0<\lambda_i<1$ for
$i=1,\ldots,d$. Thus we get $\alpha_i>0$ for $i=1,\ldots,d$, a contradiction. 
If $\dim(P)=d$, then $\alpha\in P^{\rm o} 
\subset S^{\rm o}$. As in the previous case, but now using that $S$
is a $d$-simplex, 
we get $\alpha_i>0$ for all $i$, a
contradiction. Hence $\alpha\in\partial P$. \QED  

\begin{Lemma}{\rm \cite[p.~38]{aigner}}\label{inc-exc}
Let $A_1,\ldots,A_t$ be finite subsets of a set $S$,
then
\begin{small}
$$
\left|\bigcup_{i=1}^t
A_i\right|=\sum_{i=1}^t|A_i|-\sum_{i<j}^t|A_i\cap A_j|+
\sum_{i<j<k}^t|A_i\cap A_j\cap
A_k|\mp\cdots+(-1)^{t-1}\left|\bigcap_{i=1}^t A_i\right|
$$
\end{small}
\end{Lemma}

\begin{Proposition}\label{boletin} Let $I_i$ be the ideal obtained from $I$ by 
making $x_i=0$ and let $e(I_i)$ be its multiplicity. If 
$c_dx^d+c_{d-1}x^{d-1}+\cdots+ c_1x+c_0$ is the Hilbert polynomial of the 
filtration $\mathcal{F}=\{\overline{I^n}\}_{n=0}^\infty$, then 
$$2c_{d-1}\geq
\sum_{i=1}^{d-1}\frac{e(I_i)}{(d-1)!}.
$$
\end{Proposition}

\demo  Case (I): $\dim(P)=d$. Let 
$E_{S}(x)=a_dx^d+\cdots+a_1x+1$ (resp. 
$E_{P}(x)=b_dx^d+\cdots+b_1x+1$) be the Ehrhart polynomial of $S$
(resp. $P$). By Proposition~\ref{may29-1-05}, we have 
the equality $c_i=a_i-b_i$ for all $i$. 
From the decompositions
$$
P=P^{\rm o}\cup \partial P,\ \ S=S^{\rm o}\cup \partial S, \ \ 
\partial S=K\cup H,\ \ 
\partial P= L\cup H,
$$
and using the reciprocity law of Ehrhart (Remark~\ref{ehrhart-prop}) we get:
\begin{eqnarray*}
f(n)&=&E_{S}(n)-E_{P}(n)\\
&=&E_{S}^{\rm o}(n)+
|\partial(nS)\cap\mathbb{Z}^d|-(E_{P}^{\rm o}(n)+
|\partial(nP)\cap\mathbb{Z}^d|)\\ 
&=&(-1)^dE_{S}(-n)-(-1)^dE_{P}(-n)+
|nK\cap\mathbb{Z}^d|-
|nL\cap\mathbb{Z}^d|
\end{eqnarray*}
for $0\neq n\in\mathbb{N}$. Therefore, after simplifying this
equality, we obtain: 
$$
2(c_{d-1}n^{d-1}+c_{d-3}n^{d-3}+\mbox{terms of lower degree})=
|nK\cap\mathbb{Z}^d|-|nL\cap\mathbb{Z}^d|=g(n).
$$
By the comments just before Lemma~\ref{injective}, we have the inclusions:
$$
\psi(L)\subset{M:=}\left[\left(\bigcup_{i=1}^{d-1}(\partial{P}\cap
K_i)\right)
\cup K_d\right]\setminus H \subset
K:=\left(\bigcup_{i=1}^dK_i\right)\setminus H.
$$
Using Lemma~\ref{may26-05}, we obtain:
$$
M=\left(\bigcup_{i=1}^{d-1}({P}_i\setminus 
H_i)\right)
\cup (K_d\setminus H_d)\ \mbox{ and }\ 
K=\bigcup_{i=1}^d(K_i\setminus H_i).
$$
Set $h(n)=
|nK\cap\mathbb{Z}^d|-|nM\cap\mathbb{Z}^d|$. Since $P_i\setminus
H_i\subset P_i$, $K_i\setminus
H_i\subset K_i$  for all $i$ and because $P_i\cap P_j$, $K_i\cap K_j$
are polytopes of dimension at most $d-2$ for $i\neq j$, 
by the inclusion-exclusion principle (Lemma~\ref{inc-exc})
we obtain:
\begin{eqnarray*}
&&h(n)=\sum_{i=1}^{d}|n(K_i\setminus H_i)\cap\mathbb{Z}^d|-
\sum_{i=1}^{d-1}|n(P_i\setminus H_i)\cap\mathbb{Z}^d|\\ 
& &\ \ \ \ -|n(K_d\setminus
H_d)\cap\mathbb{Z}^d|+p(n)=
\sum_{i=1}^{d-1}(E_{K_i}(n)-E_{P_i}(n))+p(n)\ \ \ \ (n\gg 0),
\end{eqnarray*}
where $|p(n)|$ is bounded by a polynomial function $P(n)$ of degree
at most $d-2$. In particular
$\lim_{n\rightarrow\infty}(p(n)/n^{d-1})=0$. By
Lemma~\ref{injective}, the 
map $\overline{\psi}\colon nL\rightarrow n\psi(L)$ given by
$\overline{\psi}(n\alpha)=n\psi(\alpha)$ is injective. Hence 
$$  
\overline{\psi}(nL\cap\mathbb{Z}^d)\subset n\psi(L)\cap
\mathbb{Z}^d\subset nM\cap\mathbb{Z}^d\ \Rightarrow\  
|nL\cap\mathbb{Z}^d|\leq |nM\cap\mathbb{Z}^d|.
$$
Consequently $g(n)=|nK\cap\mathbb{Z}^d|-
|nL\cap\mathbb{Z}^d|\geq h(n)$. Altogether we get 
$$
2c_{d-1}=\lim_{n\rightarrow\infty}\frac{g(n)}{n^{d-1}}\geq
\lim_{n\rightarrow\infty}\frac{h(n)}{n^{d-1}}=  
\lim_{n\rightarrow\infty}\left(\frac{\sum_{i=1}^{d-1}(E_{K_i}(n)-
E_{P_i}(n))}{n^{d-1}}+\frac{p(n)}{n^{d-1}}\right).
$$ 
Therefore the required inequality follows by observing that the 
polynomial function $f_i(n)=E_{K_i}(n)-E_{P_i}(n)$ is the Hilbert
function of $I_i$. Thus $f_i(n)$ has degree $d-1$ and its 
leading coefficient is equal to $e(I_i)/(d-1)!$.

Case (II): $\dim(P)=d-1$. Let $E_{S}(x)=a_dx^d+\cdots+a_1x+1$ (resp. 
$E_{P}(x)=b_{d-1}x^{d-1}+\cdots+b_1x+1$) be the Ehrhart polynomial of $S$
(resp. $P$). There is an injective map from $nP$ 
to $nK_d$ induced by $\alpha\mapsto(\alpha_1,\ldots,\alpha_{d-1},0)$.
Hence 
$$
{\rm vol}(P)=
\lim_{n\rightarrow\infty}\frac{|\mathbb{Z}^d\cap nP|}{n^{d-1}}
\leq \lim_{n\rightarrow\infty}\frac{|\mathbb{Z}^d\cap nK_d|}{n^{d-1}}=
{\rm vol}(K_d).
$$
The facets of $S$ are $K_1,\ldots,K_d$ and $K_{d+1}:=P$. Therefore 
by Proposition~\ref{may29-1-05} and using the formulas for $a_{d-1}$
and $b_{d-1}$ (see Remark~\ref{ehrhart-prop}) we conclude:
\begin{eqnarray*}
c_{d-1}=a_{d-1}-b_{d-1}&=&\frac{1}{2}\sum_{i=1}^{d+1}{\rm
vol}(K_i)-{\rm vol}(P) 
=-\frac{1}{2}{\rm vol}(P)+\frac{1}{2}\sum_{i=1}^d{\rm vol}(K_i)\\ 
&\geq& 
\frac{1}{2}\sum_{i=1}^{d-1}{\rm
vol}(K_i)=\frac{1}{2}\sum_{i=1}^{d-1}\frac{e(I_i)}{(d-1)!}.
\end{eqnarray*}\QED 

\medskip

Let $e_0,e_1,\ldots,e_d$ be the Hilbert coefficients of $f$. Recall
that we have:
\begin{small}
$$
f(n)=e_0{n+d-1\choose d}-e_1{n+d-2\choose
d-1}+\cdots+(-1)^{d-1}e_{d-1}{n\choose 1}+(-1)^de_d,
$$
\end{small}
where $e_0=e(I)$ is the multiplicity of $I$ and $c_d=e_0/d!$. Notice
that $e_d=0$ because $f(0)=0$, and $e_i\geq 0$ for all $i$, this
follows from \cite{marley}.

\begin{Corollary}\label{dec10-06} $e_0(d-1)-2e_1\geq
e(I_1)+\cdots+e(I_{d-1})\geq d-1$. 
\end{Corollary}

\demo  From the equality $c_{d-1}=\frac{1}{d!}\left[e_0{d\choose 2}-de_1\right]$
and using Proposition~\ref{boletin} we obtain the desired 
inequality. \QED 

\begin{Example}\rm Let $\mathfrak{m}=(x_1,\ldots,x_d)$ and let 
$I=\mathfrak{m}^k$. Then 
\begin{small}
$$
f(n)={kn+d-1\choose
d}=\frac{k^{d}}{d!}n^d+\frac{k^{d-1}}{(d-2)!2}n^{d-1}+
\mbox{ terms of lower degree},
$$
\end{small}
$e_0=k^d$, 
$e_1=(d-1)(k^d-k^{d-1})/2$, and we have equality in Proposition~\ref{boletin}. 
\end{Example}

\medskip

\begin{center}
ACKNOWLEDGMENT
\end{center}

\noindent The author thanks Wolmer Vasconcelos for many stimulating 
discussions.

\bibliographystyle{amsplain}

\end{document}